\documentclass[11pt,a4paper]{article}
\usepackage{wrapfig,epsfig,latexsym,amssymb}

\newcommand{\R}{\mathbb{R}}

\long\def\symbolfootnote[#1]#2{\begingroup\def\thefootnote{\fnsymbol{footnote}}
\footnote[#1]{#2}\endgroup}

\begin{document}

\thispagestyle{empty}

\title{A Very Easy Approach to the Holomorphic Discrete Series}

\author{Adam Kor\'{a}nyi}

\date{}

\maketitle

\abstract{Based on the material in an old article of Wallach a short proof of the Harish-Chandra condition is given.}

\bigskip
\bigskip


{\bf Introduction.}  Let $G$ be a simply connected simple Lie group of Hermitian type and $\Omega \simeq K\backslash G$ 
the corresponding bounded symmetric domain. The holomorphic discrete series of representations of $G$ is defined
in [HCV] and [HCVI]. As it became clear later, it is the restriction to the (uniquely meaningful)
holomorphic subspace $H^{\pi}$ of the unitarily induced representation starting from (any fixed)
irreducible unitary representation $\pi$ of $K$. The difficult question is whether $H^{\pi}$ is
non-zero: It is answered by Harish-Chandra's criterion in [HCVI].

The rather complicated original proof was, in the special case of the symplectic group, replaced by a
simpler one by Godement [G]. His argument was extended to the general case (and beyond) by Rossi and Vergne
[RV]; the method uses the Euclidean Fourier transform and the realization of $\Omega$ as a generalized halfplane.

In the present note we give a still simpler proof, staying all the time within the bounded domain
$\Omega$. (In the special case where $\pi$ is one-dimensional this proof reduces to almost nothing; see 
e.g. [FK] where it is obvious that the defining integral (3.7) of $L^{2}_{\lambda}$ converges if and
only if $\lambda > p-1$.)

A very large part of what we use is already in the article [W] of Wallach. The first part of [W] 
contains an exposition, with proofs, of the construction of the holomorphic discrete series, incorporating
a lot of previous work. We take [W] as our basic reference; we don't repeat its notations and lemmas,
just refer to them by number.

Before beginning our arguments, we point out two things: First, instead of the reference to Satake in [W]
for the critical property of $\tilde{\mathcal{K}}$ and $\tilde{D}$ it
may be easier to use the book [S], pp. 64--66. Second, instead of the non-trivial representation theoretic proof
that if $H^\pi$ is non-zero then it contains all constant functions (Lemmas 2.6, 2.7 in [W]) one can give an
alternative ``analytic" proof by noting that $H^\pi$ has a reproducing kernel which, by group invariance
properties is equal to $c({\pi}(D(z_1:z_2))^{-1}$, (with some $c>0$), and, since $D(z:0)\equiv e$, this kernel
imbeds the constants into $H^\pi$.

{\bf The Harish-Chandra condition.} We use the definitions and results of [W] with one substantive change: we
take $\gamma_r$  to be the 
highest root of $\mathfrak g$ and construct $\gamma_{r-1},\gamma_{r-2},... $ by reversing
the original process.

We also keep the notations of [W] except for abbreviating $H_{\gamma_i}$ to $h_i$, 
$e^{\lambda\Lambda_1} \otimes \pi_0$ to $\pi$, and $H^{\pi_0 \lambda}$ to $H^\pi$.

The references to numbered statements are to those in [W], unless stated otherwise.

{\bf Remark 1.} From the fact that the restricted
roots form a root system it follows easily that all multiplicities $r_{ij}$ in [W] are the same number,
to be denoted $a$, and all $r_i$ are equal to another number $b$. Also,
the $\gamma_i$ and $\alpha_1$ have the same length; they are long roots
in case there are different root lengths (cf, [H] Ch.5, \S{4}).
It is easy to check that
$$
2\varrho_P(h_i) = p \hskip 0.8cm (\forall i)
$$
and
$$
\varrho(h_r) = p-1
$$
with $p=(r-1)a+b+2$.

As it is proved in the Appendix of [RV] by a direct elementary argument, the original Harish-Chandra
condition is equivalent to the condition in the following theorem.

{\bf Theorem.} Let $\pi=e^{\lambda\Lambda_1} \otimes \pi_0$ be an 
irreducible unitary representation of $K$, $\Lambda_0$
the highest weight of $\pi_0$. Then $H^\pi$ is non-zero if and only if
$$
\lambda<-(\Lambda_0+\varrho)(h_r).
$$

{$\mathit Proof.$} By Cor.~2.8, $H^\pi = H^{\pi_0 \lambda}$ is non-zero if and only if the integral
\begin{equation}
\int_\Omega(e^{\lambda\Lambda_1} \otimes \pi_0)(D(z:z))\mathrm{d}\mu(z)
\end{equation}
converges. Since the integrand is a positive definite
$Hom(V,V)$-valued function, the integral is convergent if and only if the
integral of the $V$-trace is convergent. We use ``polar coordinates''
$z = \mathrm{Ad}(k)\sum_it_iX_i$ (as in the proof of Lemma 3.1).
In terms of these, by Lemma 2.1, we have the integral of the $V$-trace of
\begin{equation}
K(z) = (e^{\lambda\Lambda_1} \otimes \pi_0)(\exp -\sum_i\log(1-t_i^2)h_i).
\end{equation}
(Note that $K(z)$ is independent of $k$.)

We choose an orthonormal basis $\{v^s\}$ in $V$ such that each $v^s$
is a weight vector of $\pi_0$. We denote the corresponding weight by
$\Lambda^s$. Then
\begin{equation}
\mathrm{tr}_V K(z) = \sum_s \left< K(z)v^s,v^s \right>.
\end{equation}
This is a sum of positive terms; its integral is convergent (i.e. finite) if and only if
it is convergent for each term, i.e. 
\begin{equation}
\int_\Omega \left< K(z)v^s,v^s \right> \mathrm{d}\mu(z) < +\infty \hskip 0.8cm(\forall s).
\end{equation}
(Note that (4) is just the integral ($\star\star$) in the proof of Lemma 3.1 with
 $f(z) \equiv v^s$
and with $\pi_0(k)$ left out.) The integral formula used in the proof of Lemma 3.1 now
gives for the left hand side of (4):
$$
\int_{K_0}\mathrm{d}k_0 \int\mathrm{d}t_1...\mathrm{d}t_rP(t_1,...,t_r) \prod_i(1-t_i^2)^{-(\lambda\Lambda_1+\Lambda^s+2\varrho_P)(h_i)}
$$
with a polynomial $P$.
This integral is finite if and only if each of the exponents is $>$$-1$. Using Remark 1 and the fact that
$\Lambda_1(h_i)=1$ (by definition of $\Lambda_1$ and since $\gamma_i$ and $\alpha_1$ have the
same length), we get that $H^{\pi}\not=0$ if and only if
\begin{equation}
\lambda<1-p-\Lambda^s(h_i) \,\,\,\,(\forall s,\forall i)
\end{equation}

It remains to prove only that (5) is equivalent to the single condition
\begin{equation}
\lambda<1-p-\Lambda_0(h_r).
\end{equation}
In fact, by Remark 1, (6) is equivalent to the statement of our theorem.

To do this, we have to prove that $\Lambda^s(h_i)\le\Lambda_0(h_r)$, or, what
amounts to the same (since each $\gamma_i$ has the same length),
\begin{equation}
\left< \Lambda^s,\gamma_i \right> \le \left< \Lambda_0,\gamma_r \right> \hskip 0.8cm (\forall s,\forall i).
\end{equation}

This is in the real Euclidean space $i\mathfrak h_{\star}$ which is an orthogonal sum
$\R H_1 \oplus i {\mathfrak h}'_\star$ with ${\mathfrak h}'_\star$ maximal abelian
in $\mathfrak k_1$. We decompose both sides of (7) following this direct sum.

All $\Lambda^s$ are the same on $H_1$ (which is in the center of ${\mathfrak k}_C$)
by irreducibility of $\pi_0$. Also, the restriction of $\gamma_i$ to $\R H_1$ is the same for all
$i$, since $\gamma(H_1)=1$ for all $\gamma\in\Delta_P^+$. So we need only prove
\begin{equation}
\left< \bar{\Lambda}^s,\bar{\gamma}_i \right> \le \left< \bar{\Lambda}_0,\bar{\gamma}_r \right>
\end{equation}
where the bars denote restriction to $i\mathfrak h'_{\star}$. Now the $\bar{\gamma}_i$ 
are weights of the irreducible representation $\mathrm{Ad}_{\mathfrak p^+}$ of the
derived group $K_1$, with $\bar{\gamma}_r$ the highest weight. So (8) is a special case
of a general inequality for weights of any two irreducible representations: Using a Weyl group element on the left 
we may assume that $\bar{\Lambda}^s$ is dominant, in which case (8) is trivial. \hskip 0.2cm $\Box$

\end{document}